\documentclass{article} 
\usepackage{graphicx}

\begin{document}

\newtheorem{thm}{Theorem}
\newtheorem{conj}{Conjecture}
\newtheorem{lem}{Lemma}
\newtheorem{prop}{Proposition}
\newtheorem{cor}{Corollary}

\title{ Simple Closed Geodesics in Hyperbolic 3-Manifolds}
\author{Colin Adams, Joel Hass and Peter Scott\footnote{
This paper was completed while the authors were visiting the Mathematical
Sciences Research Institute in Berkeley in 1996/7. Research at MSRI is
supported in part by NSF grant DMS-9022140. The third author is grateful for
the partial support provided by NSF grants DMS-9306240 and DMS-9626537.}}
\maketitle

\begin{abstract}
This paper determines which orientable hyperbolic 3-manifolds
contain simple closed geodesics. The Fuchsian group corresponding
to the thrice-punctured sphere generates the only example of a complete
non-elementary orientable hyperbolic 3-manifold that does not contain a
simple closed geodesic. We do not assume that the manifold is geometrically
finite or that it has finitely generated fundamental group.
\end{abstract}

%AMS Classification # 53C22


\section{Introduction}
The question of which Riemannian manifolds admit simple closed geodesics is
still a mystery. It is not known whether all closed Riemannian manifolds
contain simple closed geodesics. For closed manifolds with nontrivial
fundamental group, a simple closed geodesic can always be found by taking
the shortest homotopically nontrivial closed geodesic. When the manifold is
closed but simply connected, the question is open for dimensions three and
above. In dimension two, it is known by the theorem of Lusternik and
Schnirelmann \cite{LS}, see also \cite{Grayson} and \cite{HassScott} that
the 2-sphere equipped with any smooth Riemannian metric contains at least
three distinct simple closed geodesics.

Non-compact manifolds do not necessarily contain closed geodesics, Euclidean
space being an obvious example. Even if the manifold is not simply
connected, it may not contain any simple closed geodesics, as with the
hyperbolic thrice-punctured sphere. However among the orientable, finite
area, complete hyperbolic 2-manifolds, the thrice-punctured sphere is the
only example that contains no simple closed geodesic. In this paper, we will
determine which orientable hyperbolic 3-manifolds do and do not contain
simple closed geodesics. We will prove that the Fuchsian group corresponding
to the thrice-punctured sphere generates the only example of a complete
non-elementary orientable hyperbolic 3-manifold that does not contain a
simple closed geodesic. We do not assume that the manifold is geometrically
finite or even that it has finitely generated fundamental group. The simple
closed geodesic which we produce arises from an interesting class of
elements of the fundamental group. It is the shortest closed geodesic
corresponding to a screw motion induced by the action of the fundamental
group on hyperbolic 3-space.

In proving our result, we use geometric methods to obtain results about
isometries of hyperbolic 3-space. These results apply to elliptic, as well
as to parabolic and loxodromic isometries, and we state them in full
generality.

A related question is whether a hyperbolic 3-manifold always contains a
non-simple closed geodesic. In \cite{CR}, Alan Reid and Ted Chinburg
utilized arithmetic hyperbolic 3-manifold theory to construct examples of
closed hyperbolic 3-manifolds in which every closed geodesic is simple.

We will be working with orientable 3-manifolds, so all of the isometries
discussed will be orientation preserving. A {\em screw motion} in hyperbolic
space is an isometry that has two fixed points on ${S_{\infty }^{2}}$, and
is the composition of a nontrivial translation and nontrivial rotation along
the geodesic with those fixed points as endpoints. In the language of
Kleinian groups, it is a loxodromic isometry that is not hyperbolic. A {\em
non-screw motion} will be any isometry of hyperbolic space that is not a
screw motion. Note that an isometry is a non-screw motion if and only if it
preserves an oriented plane in $H^{3}$. Moreover, if an orientation
preserving isometry is represented by a matrix in $SL(2,C),$ then it is a
non-screw motion if and only if its trace is real.

A natural place to look for a simple closed geodesic is to take the shortest
essential closed curve in a manifold, if such a curve exists. If it has
self-intersections, one can perform a cut and paste at a singular point to
create shorter curves. The problem with this line of argument is that these
shorter curves may represent parabolic elements of the fundamental group,
and thus not be homotopic to a closed geodesic. The basic method of our
argument is to show that the shortest screw motion, if it exists, does in
fact give a simple closed geodesic. In geometrically infinite manifolds, a
shortest screw motion may not exist, but an approximation argument extends
the proof to this case.

\section{Non-screw motions of $H^3$}

In this section we examine the subgroups of the isometry group of $H^3$
which are generated by elements that preserve some oriented plane.

\begin{lem}
\label{non-screw} Let $\alpha $, $\beta $ and $\alpha \beta $ be non-screw
motions of $H^{3}$. If $\alpha $ and $\beta $ do not share a fixed point in $
H^{3}\cup {S_{\infty }^{2}}$, then they preserve a common oriented plane,
and the group $<\alpha ,\beta >$ consists of non-screw motions which
preserve this plane.
\end{lem}

{\bf Proof:} We will give a geometric proof of this result. Every
orientation preserving isometry of $H^{3}$ can be represented as the product
of $180^{o}$ rotations about two axes in hyperbolic space. Call such a
rotation a {\em half-turn.} If the isometry is parabolic, both half-turn
axes pass through the parabolic fixed point. The second endpoint of one of
the half-turn axes can be chosen arbitrarily. If the isometry is hyperbolic
or elliptic, both half-turn axes pass perpendicularly through the axis of
the isometry. One of the half-turn axes can be chosen arbitrarily subject to
this property. If the isometry is elliptic, the two half-turn axes intersect
the axis of the given isometry in the same point. Let $r_{l}$ denote the
half-turn with axis $l.$ Our hypothesis that $\alpha $ and $\beta $ have no
common fixed points implies that, in any combination of these cases, we can
choose half-turns $r_{\alpha _{1}},r_{\alpha _{2}}$ for $\alpha $ and $
r_{\beta _{1}},r_{\beta _{2}}$ for $\beta $ such that $r_{\alpha
_{2}}=r_{\beta _{1}}$. This choice determines the half-turns $r_{\alpha
_{1}} $ and $r_{\beta _{2}}$. Then $\alpha \beta =r_{\alpha _{1}}r_{\alpha
_{2}}r_{\beta _{1}}r_{\beta _{2}}=r_{\alpha _{1}}r_{\beta _{2}}$. Note that
if a hyperbolic isometry $\alpha $ is expressed as a product of two
half-turns $r_{\alpha _{1}},r_{\alpha _{2}}$ as above, then $\alpha $ is
non-screw if and only if $\alpha _{1}$ and $\alpha _{2}$ lie in a common
plane.

Suppose that all three of $\alpha $, $\beta $ and $\alpha \beta $ are either
parabolic or elliptic. In order that $\alpha \beta $ be parabolic or
elliptic, it must be that $\alpha _{1}$ and $\beta _{2}$ have a point in
common, possibly on the sphere at infinity. Thus, the three axes $\alpha
_{1},\alpha _{2}$ and $\beta _{2}$ create a triangle in $H^{3}$ and the
product of any two of the half-turns about them preserves the oriented plane
containing that triangle.

Suppose next that exactly two of $\alpha $, $\beta $ and $\alpha \beta $ are
parabolic or elliptic. Since the inverse of a parabolic or elliptic element
is parabolic or elliptic, we can rename the elements so that $\alpha $ and $
\beta $ are parabolic or elliptic. Then $\alpha _{1}$ and $\alpha _{2}$
share a common point, as do $\beta _{1}(=\alpha _{2})$ and $\beta _{2}$.
Since $\alpha \beta $ is non-screw, $\alpha _{1}$ and $\beta _{2}$ lie in a
common plane. As two distinct points on $\alpha _{2}$ also lie in this
plane, so does $\alpha _{2}$ itself. As before, the product of any two of
the half-turns about $\alpha _{1},\alpha _{2}$ or $\beta _{2}$ preserve this
oriented plane.

In the case where exactly one of $\alpha $, $\beta $, $\alpha \beta $ is
parabolic or elliptic, we can assume that it is $\alpha $. Then $\alpha _{1}$
and $\alpha _{2}$ share a common point, while $\beta _{1}(=\alpha _{2})$ and 
$\beta _{2}$ lie in a common plane $P.$ Moreover $\alpha _{1}$ and $\beta
_{2}$ also lie in a common plane $Q.$ If $P$ and $Q$ are distinct, they
intersect in a line passing through $\alpha _{1}\cap \alpha _{2},$ and this
line must equal $\beta _{2}.$ Thus $\alpha _{1},\alpha _{2}$ and $\beta _{2}$
have a common point, which must be fixed by both $\alpha $ and $\beta ,$
contradicting our hypothesis. If $P=Q,$ then $\alpha _{1},\alpha _{2}$ and $
\beta _{2}$ all lie in this plane, and the product of any two of the
half-turns about them will preserve this oriented plane as required.

The final case occurs when none of $\alpha $, $\beta $ and $\alpha \beta $
are parabolic or elliptic. In this case each pair of axes from $\alpha
_{1},\alpha _{2}$ and $\beta _{2}$ are disjoint, even at infinity, but share
a common plane. Let $P$ denote the plane containing $\alpha _{1}$ and $
\alpha _{2},$ let $Q$ denote the plane containing $\alpha _{2}$ and $\beta
_{2},$ and let $R$ denote the plane containing $\alpha _{1}$ and $\beta
_{2}. $ Let $\lambda $ denote the common perpendicular line segment between $
\alpha _{1}$ and $\alpha _{2}.$ Let $\Pi $ denote the plane which contains $
\lambda $ and is orthogonal to $P.$ Then $\Pi $ is also orthogonal to each
of $Q$ and $R.$ If $Q\neq R,$ it follows that $\Pi $ is orthogonal to $Q\cap
R=\beta _{2}.$ The product of any two of the half-turns about $\alpha
_{1},\alpha _{2}$ and $\beta _{2}$ will preserve this oriented plane as
required. If $Q=R,$ then $\alpha _{1},\alpha _{2}$ and $\beta _{2}$ all lie
in this plane, and as before the product of any two of the half-turns about
them will preserve this oriented plane. This completes the proof of Lemma 1.

\begin{lem}
If $\Gamma $ is a non-elementary Kleinian group and if every element of $
\Gamma $ is a non-screw motion, then $\Gamma $ is Fuchsian. \label{Fuchsian}
\end{lem}

{\bf Proof:} A proof appears in \cite{M}, [p.108]. We give here a geometric
proof. Since $\Gamma $ is not elementary, there exist infinitely many
distinct axes for hyperbolic (i.e. non-screw loxodromic) elements. Choose
two hyperbolic elements with distinct axes. Then Lemma~\ref{non-screw} shows
that the two axes lie in a plane $P$ that is preserved by the subgroup that
they generate. Suppose that there exists a limit point of $\Gamma $ that is
not contained in the circle that bounds $P.$ There is a hyperbolic axis with
an endpoint that is arbitrarily close to this point. Let $\gamma $ be a
hyperbolic element with axis $A$ that has an endpoint off the circle. The
set of all geodesics on $P$ that are coplanar with $A$ form a foliation of
the plane $P$. However, there exists a hyperbolic axis in $P$ that is not
one of the geodesics in this foliation, since we can take any two of the
axes in the foliation and let the hyperbolic element of the first act a
number of times on the second. The second axis will be carried close to one
endpoint of the first, so it cannot be a geodesic of the foliation. It
follows that the whole limit set is in the boundary of $P$ and $\Gamma $ is
Fuchsian.

\section{The Main Result}

In this section we prove the following result.

\begin{thm}
\label{geom.inf.man} Let $M$ be an orientable hyperbolic 3-manifold. Then
exactly one of the following three cases occurs:

\begin{enumerate}
\item  There exists a simple closed geodesic in $M$.

\item  $M$ is the quotient of $H^{3}$ by a Fuchsian group corresponding to a
thrice-punctured sphere.

\item  The fundamental group of $M$ is elementary with zero or one limit
point.
\end{enumerate}
\end{thm}

{\bf Proof:} Let $\Gamma $ be a discrete group of isometries such that $
M=H^{3}/\Gamma $. Note that if $\Gamma $ is elementary with two limit
points, then there exists a hyperbolic isometry such that the endpoints of
its axis are these two limit points. In order that these be the only limit
points, it must be that this is the only axis of an isometry of $\Gamma $ in 
$H^{3}$. This means that this axis projects to a simple closed geodesic. We
now assume that $\Gamma $ is non-elementary. If $\Gamma $ is Fuchsian, then
either $\Gamma $ is the fundamental group of a thrice-punctured sphere, or
there is a totally geodesic surface in $M$ that contains a simple closed
geodesic. Hence, the theorem holds in this case. We now assume that $\Gamma $
is a non-elementary non-Fuchsian discrete group. By Lemma~\ref{Fuchsian} we
can find an element in $\Gamma $ that is a screw motion. If $\gamma $ is a
screw motion we let $l(\gamma )$ denote the translation length of $\gamma ,$
and define $I=\inf \{l(\gamma ):\gamma \mbox{
is a screw motion in }\pi _{1}(M)\}$. There are now two cases depending on
whether or not $I$ is realized as the translation length of a screw motion.

{\bf Case 1:} There is a screw motion $\gamma $ with $l(\gamma )=I.$ We will
denote the corresponding closed geodesic by $c.$ We claim that $c$ must be
embedded. If $c$ is not embedded, then $c$ intersects itself transversely at
some point. We can cut and paste at this point to produce curves $a,b$ and $
ab,$ each shorter than $c,$ with $c$ homotopic to $ab^{-1}$. Then $a$, $b$
and $ab$ must correspond to non-screw motions, as otherwise we could find
closed geodesics homotopic to them that are shorter than $c$, contradicting
the minimizing property of $c$. Lemma~\ref{non-screw} now implies that the
isometry $\gamma $ corresponding to $c$ must also be a non-screw motion.
This contradiction shows that $c$ must be embedded as claimed, so that $M$
admits a simple closed geodesic.

{\bf Case 2: }There is no screw motion $\gamma $ with $l(\gamma )=I.$

Note that in the manifolds which we are considering the fundamental group
may not be finitely generated, and there need not exist a shortest closed
geodesic nor a shortest screw motion. Explicit examples of hyperbolic
3-manifolds containing no shortest closed geodesic have been constructed by
Bonahon and Otal \cite{B-O}.

In the following, we suppose that $M$ does not contain a simple closed
geodesic and will obtain a contradiction. Let ${c_{n}}$ be a sequence of
closed geodesics corresponding to screw motions whose lengths approach $I$
as $n\rightarrow \infty $. Each of these geodesics must itself be singular.
Therefore each geodesic can be subdivided at a singular point into two
sub-loops, which we call $a_{n}$ and $b_{n}$. Let $\alpha _{n}$, $\beta _{n}$
and $\gamma _{n}$ be the closed geodesics homotopic to $a_{n}$, $b_{n}$ and $
a_{n}b_{n}^{-1}$, if they exist. Let $\epsilon $ denote the Margulis
constant in dimension three, so that any closed geodesic of length less than 
$\epsilon $ must be embedded. It follows that $I\geq \epsilon $. We also
claim that there exists a number $K>0$ such that $l(a_{n})>K$ and $
l(b_{n})>K $ for all $n$. If not, then by passing to a subsequence we can
assume that $\displaystyle\lim_{n\rightarrow \infty }l(a_{n})=0$. For any
value of $n$ with $l(a_{n})<\epsilon $, the isometry corresponding to $a_{n}$
must be parabolic. However, a very short parabolic loop must occur far out
in a cusp, forcing $c_{n}$ to be very long since it is not contained
completely in a cusp. More explicitly, a shortest curve in a maximal
horosurface (which must be an annulus or torus) cutting off a cusp has
length at least 1. The shortest nontrivial curve that comes to within a
distance $d$ of the horosurface can be computed by standard hyperbolic
geometry. Inverting to get the distance in terms of the length of the curve
shows that a homotopically non-trivial curve in a cusp whose length is less
than $\delta <1$ has distance at least $\ln (\mbox{csch}(\delta /2)/2)$ from
the maximal horosurface bounding the cusp. In particular, if $a_{n}$
represents a parabolic and $\displaystyle\lim_{n\rightarrow \infty
}l(a_{n})=0$ then $\displaystyle\lim_{n\rightarrow \infty }l(c_{n})=\infty $
, contradicting our choice of ${c_{n}}$.

Now let $\theta _{n}$ be the angle between the two strands of $b_{n}$
meeting at the singular point on $c_{n}$. We will show that $\displaystyle
\lim_{n\rightarrow \infty }\theta _{n}=0$. Suppose that $\theta _{n}>\theta
_{0}>0$ for all $n$. Since $l(a_{n})$ and $l(b_{n})$ are both greater than $
K $, we can choose a segment of $a_{n}$ of length greater than $K,$ a
segment of $b_{n}$ of length greater than $K,$ each with one end at the
singular point of $c_{n}$ and consider the triangle determined by their
endpoints. This triangle has angle $\pi -\theta _{n}$ at the singular point
of $c_{n}.$ The assumption that $\theta _{n}>\theta _{0}>0$ for all $n$
implies that the third edge of this triangle is shorter than the
corresponding path along $a_{n}$ and $b_{n}$ by a constant $L>0$, depending
only on $K$ and $\theta _{0}.$ In particular, $L$ is independent of $n$. If
there is a closed geodesic $\gamma _{n}$ homotopic to the loop $
a_{n}b_{n}^{-1},$ it follows that $l(\gamma _{n})<l(c_{n})-L,$ as $
a_{n}b_{n}^{-1}$ and the original loop $c_{n}$ have the same length. If $
l(c_{n})<I+L$, it follows that either $a_{n}b_{n}^{-1}$ represents a
parabolic, or that $l(\gamma _{n})<I.$ In either case, $a_{n}b_{n}^{-1}$
represents a non-screw motion. Now if $l(c_{n})<I+K,$ the fact that $
l(a_{n})>K$ and $l(b_{n})>K$ for all $n$ implies that $l(a_{n})$ and $
l(b_{n})$ are each less than $I,$ so that $a_{n} $ and $b_{n}$ each
represent non-screw motions$.$ Thus if $n$ is large enough to ensure that $
l(c_{n})<I+L$ and $l(c_{n})<I+K,$ then $a_{n},b_{n}$ and $a_{n}b_{n}^{-1}$
will all represent non-screw motions. But now Lemma~\ref{non-screw} implies
that $c_{n}$ represents a non-screw motion -- a contradiction.

Finally, consider the case where the angle $\theta _{n}$ at the singular
point of the closed geodesic $c_{n}$ satisfies $\displaystyle
\lim_{n\rightarrow \infty }\theta _{n}=0$. Recall that there exists a number 
$K>0$ such that $l(a_{n})>K.$ As the sequence of closed geodesics $c_{n}$
has length bounded above and it is trivial that $l(a_{n})<l(c_{n})$, it
follows that we can bound the length of $a_{n}$ between two constants $
0<K<l(a_{n})<M$. Lifting to the universal cover, we have two geodesic lifts
of $c_{n}$ crossing at a pre-image $Q$ of the singular point. Along each
geodesic at a distance $l(a_{n})<M$ from the crossing point, we have two
points which also project down to the singular point. If $P$ is one of these
points, then ${a_{n}}^{2}P$ is the other, where ${a_{n}}$ denotes the
element of $\Gamma $ which sends $P$ to $Q.$ Recall that we are assuming
that as $n$ approaches $\infty $, $\theta _{n}$ approaches 0. Thus $
\displaystyle\lim_{n\rightarrow \infty }d(P,{a_{n}}^{2}P)=0$. Since
sufficiently short closed geodesics are simple, it must be that $a{_{n}}^{2}$
, and therefore $a_{n}$ itself, are parabolic for large $n$. For a parabolic
element, $d(P,a{_{n}}^{2}P)\rightarrow 0\Rightarrow d(P,a_{n}P)\rightarrow 0$
, which contradicts the fact that $l(a_{n})>K$. This contradiction completes
the proof that $M$ must contain a simple closed geodesic. Note that, in
fact, our argument shows more than this. It shows that all sufficiently
short closed geodesics in $M$ coming from screw motions must be simple.

\begin{flushleft}
Colin Adams, Department of Mathematics, Williams College,
Williamstown, MA,01267. e-mail: Colin.C.Adams@williams.edu \\
Joel Hass, Department of Mathematics, University of California,
Davis, CA 95616. e-mail: hass@math.ucdavis.edu\\
Peter Scott, Department of Mathematics, University of Michigan,
Ann Arbor MI, 48109. e-mail: pscott@math.lsa.umich.edu\\
\end{flushleft}


\begin{thebibliography}{99}
\bibitem{B-O}  F. Bonahon and J. P. Otal, Varietes hyperboliques a
geodesiques arbitrairement courtes. Bull. London Math. Soc. 20 (1988), no.
3, 255--261.

\bibitem{CR}  T. Chinburg and A. Reid, Closed hyperbolic $3$-manifolds whose
closed geodesics all are simple. J. Differential Geom. 38 (1993), 545--558.

\bibitem{Grayson}  M. A.{\ {Grayson}, {Shortening embedded curves}, {Ann. of
Math. (2)129}, ({1989)}, {71--111}.}

\bibitem{HassScott}  J. {Hass and P. Scott}, {Shortening curves on surfaces}
, {Topology,} {33} ({1994)}, {25--43}.

\bibitem{K}  S. Kojima, Isometry transformations of hyperbolic 3-manifolds,
Topology and its Appl., 29 (1988) 297-307.

\bibitem{LS}  L. Lusternik and L. Schnirelmann, M\'{e}thodes Topologiques
dans les Probl\`{e}mes Variationnelles, Hermann and Co., Paris 1934.

\bibitem{M}  B. Maskit, Kleinian Groups, Grundlehren der Mathematischen
Wissenschaften 287, Springer-Verlag, Berlin-New York, 1988.

\bibitem{My}  R. Myers, Simple knots in compact, orientable 3-manifolds,
Trans. of A. M. S., Vol.273, No. 1,(1982),75-91.

\bibitem{S}  T. Sakai, Geodesic knots in a hyperbolic 3-manifold, Kobe J.
Math., 8 (1991) 81-87.

\bibitem{Th}  W. Thurston, The geometry and topology of 3-manifolds,
Princeton University Lecture Notes (1978).
\end{thebibliography}
\end{document}